\documentclass[12pt]{article}
\usepackage[T1]{fontenc}
\usepackage[utf8]{inputenc}

\usepackage[francais,english]{babel}
\usepackage{amsmath} 
\usepackage{graphicx} 
\usepackage{geometry}          
\usepackage{subfig} 
\usepackage{amssymb}
\usepackage{amsthm}
\usepackage{epstopdf}
\usepackage{enumerate}
\usepackage{makeidx}
\usepackage{hyperref}
\usepackage{pinlabel}

\hypersetup{
backref=true, 
pagebackref=true,
hyperindex=true, 
breaklinks=true, 
urlcolor= blue, 
linkcolor= blue, 
bookmarks=true, 
bookmarksopen=true, 
pdftitle={Conformal dimension  on boundary of right-angled hyperbolic buildings}, 
pdfauthor={Antoine Clais}, 
pdfsubject={} 
}

\newtheorem{theo}{Theorem}

\newtheorem{déf}[theo]{Definition}

\newcommand{\R} {\ensuremath{\mathbb{R}}}

\newcommand{\h} {\ensuremath{\mathbb{H}}}
\newcommand{\F} {\ensuremath{\mathcal{F}}}

\newcommand{\borg} {\ensuremath{\partial \Gamma}}

\newcommand{\N} {\ensuremath{\mathbb{N}}}

\newcommand{\s} {\ensuremath{\mathbb{S}}}

\newcommand{\modcomb}[2] {\ensuremath{\mathrm{Mod}_#1} (#2)}

\newcommand{\modcombg}[1] {\ensuremath{\mathrm{Mod}_p} (#1,G_k)}

\newcommand{\modcombfo} {\ensuremath{\mathrm{Mod}_p} (\mathcal{F}_0,G_k)}

\newcommand{\modcombcoxgp}[1] {\ensuremath{\mathrm{mod}_p} (#1,G_k^W)}

\newcommand{\confdim}[1] {\mathrm{Confdim}(\ensuremath{#1})}

\newcommand{\dia}[1] {\mathrm{diam} \ensuremath{\:#1}}

     \newcommand{\liste} [1] {
  \begin {itemize}#1
  \end{itemize}}

  \makeindex       

\title{Conformal dimension  on boundary of right-angled hyperbolic buildings}
\author{Antoine Clais\\Technion\\Department of Mathematics\\32000 Haifa, Israel\\ 
\texttt{aclais@tx.technion.ac.il}}
\date{\today}

\begin{document}
  \maketitle
 \selectlanguage{french}
\begin{abstract}
Dans cette note, on utilise des modules combinatoires sur le bord d'un immeuble hyperbolique à angles droits pour encadrer sa dimension conforme. La borne inférieure obtenue est optimale dans le cas des immeubles Fuchsiens.
\end{abstract} 
  \selectlanguage{english}
  \begin{abstract}
  In this note, we use some combinatorial modulus on the boundary of a right-angled hyperbolic building  to control its conformal dimension. The lower bound obtained is optimal in the case of Fuchsian buildings.
\end{abstract}

\begin{flushleft}{\bf MSC2010:}  	20F67, 20E42.\end{flushleft}


  The conformal dimension is a quasi-isometric invariant of Gromov hyperbolic spaces that has been introduced by P. Pansu in \cite{PansuDimconf}. Since then it has become a major tool used to study quasi-conformal properties of boundaries of hyperbolic groups in relation with rigidity phenomenon. In particular, in \cite{BonkKleinerConfDimGromHypergrps}, it plays a key role in proving that right-angled Fuchsian buildings satisfy a Mostow type rigidity theorem. We refer to \cite{KleinerAsymptoticGeom} and \cite{HaissinskyGeomQConf} for surveys concerning this topic and to \cite{MackayTysonConfDim} for a survey concerning the conformal dimension in the  more general context  of self-similar spaces.
 

\paragraph{Recalls} Let $(Z,d)$ be a compact arcwise connected metric space. The \emph{cross-ratio} of four distinct points  $a,b,c,d \in Z$ is \[[a:b:c:d] = \frac{d(a,b)}{d(a,c)}\cdot\frac{d(c,d)}{d(b,d)}.\]
An homeomorphism $f : Z \longrightarrow Z$ is \emph{quasi-Moebius} if there exists an homeomorphism $\phi : \R ^+ \longrightarrow \R ^+$ such that for any quadruple of distinct points  $a,b,c,d \in Z$  \[[f(a):f(b):f(c):f(d)] \leq \phi( [a:b:c:d ]).\]
 Now we assume that $Z$ is   a \emph{$Q$-Ahlfors-regular}   (AR)  for $Q>1$. This means that there exists   a constant $C>1$ such that  for any $0<R\leq  \dia{Z}$ and any $R$-ball $B\subset Z$ one has  \[C^{-1}\cdot R^Q \leq \mathcal{H}_d(B) \leq C\cdot R^Q,\]
  where  $\mathcal{H}_d(\cdot)$ denotes the Hausdorff measure  of  $(Z,d)$. Notice that, under this assumption, $Q$ is equal to  the Hausdorff dimension $\dim_\mathcal{H} (Z,d)$ of $(Z,d)$.
The \emph{Ahlfors-regular conformal gauge} of $(Z,d)$ is defined as follows
\[\mathcal{J}_c(Z,d):= \{ (Z,\delta) : (Z,\delta) \text{ is AR and quasi-Moebius homeomorphic to }  (Z,d)\}. \]
  
\begin{déf}
 The \emph{Ahlfors-regular conformal dimension} \index{Conformal dimension (Ahlfors-regular)} of  $\borg$ is  
\[\confdim{Z,d} := \inf\{\dim_\mathcal{H}(Z,\delta) : (Z,\delta) \in   \mathcal{J}_c(Z,d)  \}  . \]
\end{déf}
In the rest of the  note we will simply call it the \emph{conformal dimension}. We recall that all the visual metrics on the boundary of a hyperbolic space are quasi-Moebius homeomorphic one to the other and AR. In particular, the conformal dimension is a quasi-isometric invariant of a hyperbolic space.  

 As the topological   and the Hausdorff dimensions are respectively invariant under homeomorphisms and bi-Lipschitz maps, the conformal dimension is invariant under quasi-Moebius maps. The inclusions between these three classes of maps imply the following inequalities: \[\dim_T (Z)\leq \confdim{Z,d}\leq \dim_\mathcal{H}(Z,d),\] where $\dim_T(Z)$ designate the topological dimension of $Z$.

Combining ideas of G.D Mostow, P. Pansu and M. Bourdon   with a theorem of M. Bonk and B. Kleiner (see \cite[Theorem 1.3]{BonkKleinerConfDimGromHypergrps}) it is known that if the conformal dimension of the boundary of a CAT(-1) group is realized in the conformal gauge, then the underlying CAT(-1) space  satisfies a Mostow type rigidity theorem (see \cite[Théoreme 5.11]{HaissinskyGeomQConf}).

  \paragraph{The result} The goal of this note is to prove the theorem bellow that relates the conformal dimension of the boundary of a building to the conformal dimension of the boundary of an apartment.  We refer to \cite{DavisBook}   for generalities concerning Coxeter groups and buildings.
    
  Let  $\mathcal{G}$ denote a \emph{finite simplicial graph} \emph{i.e} $\mathcal{G}^{(0)}$ is finite, each edge has two different vertices, and $\mathcal{G}$ contains no double edge. We denote  by  $\mathcal{G}^{(0)}=\{v_1, \dots, v_n\}$  the vertices of $\mathcal{G}$ and we set $S=\{s_1, \dots , s_n\}$.    If for $i\neq j$ the corresponding vertices $v_i, v_j$ are connected by an edge, then we write $v_i \sim v_j$. We denote by $W$ the the \emph{right-angled Coxeter group} whose relation graph is $\mathcal{G}$, namely 
 \[ W =   \left\langle s_i\in S \vert  s_i^{2}=1,  s_is_j=s_js_i \text{ if } v_i \sim v_j  \right\rangle .\]
For $q\geq 2$ we denote by $\Gamma_q$ the  group    defined by the following presentation\[ \Gamma_q  =   \left\langle s_i\in S \vert  s_i^{q}=1,  s_is_j=s_js_i \text{ if } v_i\sim v_j   \right\rangle. \] 
  In the rest of this note, we will assume that $\Gamma_q$ is infinite  hyperbolic with arcwise connected boundary $\borg_q$. This assumption only involves the graph $\mathcal{G}$ (see  \cite{MeierWhen} and \cite{DavisMeierTopo}). We also equip $\borg_q$ with a visual metric $d$.
  
 We recall that   $\Gamma_q$ acts by isometry properly discontinuously and cocompactly on $\Delta_q$ the building of type $(W,S)$ and of constant thickness $q$. As a consequence the boundaries $\borg_q$ and $\partial\Delta_q$ are canonically identified by a $\Gamma_q$-equivarient quasi-Moebius homeomorphism.  
We also recall that $\Delta_q$ is a CAT(-1) metric space.

 In the following, for $g\in \Gamma_q$, we designate by $\vert g \vert$ the distance, for the word metric on $\Gamma_q$ relative to the generating set $S^q$,  
 between $e$ and $g$. Then  $\tau(q) = \limsup_k \frac{1}{k} \log(\#\{g\in \Gamma_q : \vert g \vert \leq k \})$ is the growth rate of $\Gamma_q$.     In the rest of the note, we will write  $Q(q) = \confdim{\borg_q}$. Notice that, in  particular, $Q(2)=\confdim{\partial W}$.
  
\begin{theo} \label{letheoreme} There exists a constant $C>0$ independent of $q$ such that 
\[ Q(2) \cdot  \big( 1+\frac{\log(q-1)}{\tau(2)}\big) \leq Q(q) \leq C \log(q-1).   \] 
\end{theo}

  \paragraph{Consequences of the theorem}
    In general, the conformal dimension $Q(2)$ of $\partial W$ is unknown. However the topological dimension of $\partial W$  is easy to read in the graph $\mathcal{G}$ (see \cite{DavisMeierTopo}). Hence, one can use the inequality $\confdim{\partial W} \geq \dim_T(\partial W)$, to obtain an  explicit lower bound.

  The only hyperbolic buildings for which we can compute the conformal dimension of their boundaries are Fuschian buildings. In this case , the lower bound of Theorem \ref{letheoreme} is optimal (see \cite[Théorème 1.1]{BourdonImHyperDimConfRigi}). This allows us to think that this bound should be optimal for others  hyperbolic buildings.
  
  In particular, if $D$ is the regular right-angled dodecahedron in $\h^3$ and   $(W,S)$  is the right-angled Coxeter group generated by  the reflections about the faces of $D$. Then $\partial W$ is quasi-Moebius homeomorphic to the  Euclidean sphere $\s^2$ and  $Q(2)=2$. For $q\geq 3$ the corresponding building $\Delta_q$ has all chances to be Mostow rigid. Indeed, the rigidity of its apartments (that are copies of $\h^3$) should be increased by the building structure.   Moreover in \cite{ClaCLPonBRAHB} it is proved that $\partial \Delta_q$ satisfies the \emph{combinatorial Loewner property}. This property is conjecturally equivalent to admitting a metric  realizing the conformal dimension  (see \cite[Conjecture 7.5]{KleinerAsymptoticGeom}). As a consequence, it may be enough to find a visual metric of Hausdorff  dimension equal to $2+\frac{2\log(q-1)}{\tau(2)}$  on $\borg_q$  in order to prove the rigidity of this building.

Finally, an immediate computation shows that the growth rates $\tau(q)$ and $\tau(2)$ are related   by the formula   $ \tau(q)= \tau(2) +\log(q-1).$ On the other hand, growth rates of Coxeter groups can be computed by hand. For instance, if $W$ is the reflection group of $\h^3$   mentioned in the preceding paragraph $\tau(2)=\log(4+\sqrt{15})$.

 \paragraph{Combinatorial modulus}
 
 For a complete introduction on combinatorial modulus on boundaries of hyperbolic groups we refer to \cite{BourdonKleinerCLP}. We will restrict to the example given by the group $\borg_q$ equipped with a visual metric $d$.
 
For $k\geq 0$ and $\kappa > 1$, a \emph{$\kappa$-approximation of $(\borg_q,d)$ on scale $k$} \index{Approximation  on scale $k$} is a finite covering $G_k$ \index{$G_k$} by open subsets  such that for any $ v\in G_k$ there exists $z_v\in v$ satisfying the following properties:

\liste{ \item  $B(z_v,\kappa^{-1} 2^{-k}) \subset
v \subset B(z_v,\kappa 2^{-k}) $,
\item $\forall v , w \in G_k$ with $v\neq w$ one has $B(z_v,\kappa^{-1} 2^{-k})\cap B(z_w,\kappa^{-1} 2^{-k})= \emptyset$.}
A sequence $\{G_k \}_{k\geq 0}$ is called a \emph{$\kappa$-approximation of $(\borg_q,d)$}.
 
Now we fix the approximation $\{G_k \}_{k\geq 0}$. Throughout this note, we will call a curve in $\borg_q$ a continuous map $\gamma : [0,1]\longrightarrow \borg_q$ and we will designate by $\gamma$ its image. Let $\rho: G_k \longrightarrow [0,+\infty)$ be a positive function and $\gamma$ be a curve in $\borg_q$. The \emph{$\rho$-length} of   $\gamma$   is $L_\rho (\gamma) = \sum\limits_{\gamma\cap v \neq \emptyset } \rho (v).$ For $p\geq 1$, the \emph{$p$-mass} of $\rho$ is  $M_p(\rho) = \sum\limits_{v\in G_k} \rho (v)^p.$
 Let $\F$ be a non-empty set of curves in $\borg_q$. We say that the function $\rho$ is \emph{$\F$-admissible} if $L_\rho(\gamma)\geq1$ for any curve $\gamma\in \F$.  For $p\geq 1$, the \emph{$G_k$-combinatorial $p$-modulus of $\F$}   is \[ \modcombg{\F}= \inf \{ M_p(\rho)\},\] where the infimum is taken over the set of $\F$-admissible functions and with   $\modcombg{ \emptyset }:=0$.

 \paragraph{Two ingredients  for the proof} The first ingredient of the proof  is that  the conformal dimension is a critical exponent for the combinatorial modulus.  This is the immediate application of a theorem that holds in the general context of  boundaries of    hyperbolic groups (see \cite{KeithKleiner} or \cite[Theorem 0.1.]{CarrascoConfDimModComb}). 
 
  Let $d_0>0$ be a small  constant compared with the geometric constants of $(\borg_q,d)$ and let $\F_0$ be the set of all the curves in $\borg_q$ of diameter larger than $d_0$. Under this assumption\[(*) \ \ \confdim{\borg_q} = \inf \{p\in [1, +\infty) : \lim_{k\rightarrow  + \infty} \modcombfo = 0\}.\]

    The second ingredient of the proof is   a control of the combinatorial modulus on the boundary of the building by the combinatorial modulus on the boundary of an apartment (see \cite[Theorem 9.1]{ClaCLPonBRAHB}). 
  
 Let $G_k^W$ and by $G_k$ two $\kappa$-approximations of respectively $\partial W$ and $\borg$ constructed as in \cite[Subsection 8.2]{ClaCLPonBRAHB}. We  recall that by construction, there exists $\lambda\geq 1$ such that, for any $k\geq 0$ \[ \#G_k^W \leq \#\{g\in W  : \vert g\vert \leq k\} \leq  \lambda \cdot \#G_k^W.\] In particular, this implies that $\tau(2)=\limsup_{k\rightarrow\infty} \frac{1}{k} \log (\#G_k)$. Now we  designate by $\modcombcoxgp{\cdot}$ and $\modcombg{\cdot}$ the combinatorial modulus computed respectively in $\partial W$ and in $\borg_q$. Let $d_0>0$ be a small constant compared with the geometric parameters of $\partial W$ and $\borg_q$ and  let  $\F_0^W$  and  $\F_0$ be the set of all the curves   respectively in $\partial W$  and $\borg_q$  of diameter larger than $d_0$. Under these hypothesis, for any $p\geq 1$, there exists a constant $D\geq 1$ such that for every $k \geq 1$
\[(**) \ \ D^{-1}\cdot \modcomb{p}{\F_0,G_k}  \leq    (q-1)^k \cdot \modcombcoxgp{F^W_0}  \leq D \cdot  \modcomb{p}{\F_0,G_k}.\] 
 
As an immediate consequence of $(*)$ and $(**)$\[Q(q) = \inf \{p\in [1, +\infty) : \lim_{k\rightarrow  + \infty} (q-1)^k\cdot\modcombcoxgp{\F_0^W} = 0\}.\] We use this equality in proving Theorem \ref{letheoreme}. Throughout this proof, $\{\rho_k\}_{k\geq 0}$ is a sequence of $\F_0^W$-admissible functions.

  \paragraph{Right-hand side inequality}
 Comparing the diameters of the elements of $G_k$ with $d_0$ one obtains that there exist two constants  $K>0$ and $0<\lambda <1$ independent of $k$ and $q$ such that, up to change the sequence $\{\rho_k\}_{k\geq 0}$, for all $k\geq 1$  and all  $w\in G_k^W$ 
\[\rho_k(w)\leq   K \cdot \lambda^{k}.\] Hence 
\[  (q-1)^k\sum_{w\in G_k^W} \rho_k(w)^p \leq K^p\cdot   \#G_k^W \cdot [\lambda^p(q-1)]^k. \] As we recalled,    $\tau(2)=\limsup_{k\rightarrow\infty} \frac{1}{k} \log (\#G_k)$, thus $ Q(q)\leq \frac{1}{\log 1/\lambda} ( \tau(2)+ \log(q-1)).$ 
  
  \paragraph{Left-hand side inequality}
 \vspace{0.5cm}

Now we set  $G=\bigcup\limits_{k\geq 0} G_k^W$ and for  $w\in G$ we write $\vert w\vert =k$ such that  $w\in G_k^W$. For  a sequence  $\{\rho_k\}_{k\geq 0}$ of  $\F_0^W$-admissible function, let $\rho:w \in G \longrightarrow\rho_{\vert w \vert}(w)\in \R^+$.  Now we observe that \[ \text{ if } \sum_{w\in G}(q-1)^{\vert w\vert}  \rho(w)^p < \infty \text{ then } p\geq Q(q).\]
We define the function \[P_\rho(s):=\inf\{p>0 \vert  \sum_{w\in G} (q-1)^{s\vert w\vert} \rho(w)^p<\infty  \} \]
  and we study the function $P_\rho$ to compute the lower bound for $Q(q)$.
 
 To start, one has  $P_\rho(0)\geq Q(2)$ and, up to change $\{\rho_k\}_{k\geq 0}$, we can choose $P_\rho(0)$ arbitrarily close to $Q(2)$.  Then, identically, one has  $P_\rho(1)\geq Q(q)$ and, up to change $\{\rho_k\}_{k\geq 0}$, we can choose $P_\rho(1)$ arbitrarily close to $Q(q)$. 
  
  Now we set   $s_0= \sup\{s\in \R \vert \sum_{w\in G} (q-1)^{s \vert w\vert} <\infty\}$. Clearly $s_0<0$ and is such that $P_\rho(s_0)=0$. Indeed, as we said in the preceding paragraph, there exist two constants  $K>0$ and $0<\lambda <1$ independent of $k$ and $q$ such that, up to change the sequence $\{\rho_k\}_{k\geq 0}$, for all $k\geq 1$  and all  $w\in G_k^W$ 
\[\rho(w)\leq   K \cdot \lambda^{\vert w\vert}.\] Hence,  if for all  $\epsilon > 0 $ \[\sum_{w\in G} (q-1)^{-\epsilon \vert w\vert } . (q-1)^{s_0 \vert w\vert }<\infty.\]  Then for all $\epsilon'>0$ small enough one has \[\sum_{w\in G} \rho(w)^{\epsilon'} (q-1)^{s_0 \vert w\vert }\leq K^{\epsilon'} \cdot \sum_{w\in G} (q-1)^{-\epsilon \vert w\vert } . (q-1)^{s_0 \vert w\vert }<\infty. \]
On the other hand $s_0= \sup\{s\in \R \vert \sum_{k\in  \N} \#G_k^W (q-1)^{sk}<\infty\}$. As we recalled, $ \tau(2)=\limsup_k \frac{1}{k} \log  (\#G_k^W)$ thus  \[s_0=-\frac{\tau(2)}{\log(q-1)}.\]
 
Now we check that the function $P_\rho$ is convex on $[s_0,+\infty)$. In other word, we check that for all $t \in [0,1]$, for all $[a,b] \subset [s_0,+\infty)$ and for all $\epsilon >0$ one has  \[ \sum_{w\in G} (q-1)^{(t a + (1-t)b )\vert w\vert} \rho(w)^{(t P_\rho(a) + (1-t) P_\rho(b) )+\epsilon}<\infty.\]
Indeed, for all $\epsilon >0$  \[ \sum_{w\in G} x_w:=\sum_{w\in G} (q-1)^{ a \vert w\vert} \rho(w)^{ P_\rho(a)  +\epsilon}<\infty \text{ and  } \sum_{w\in G} y_w:=\sum_{w\in G} (q-1)^{ b \vert w\vert} \rho(w)^{ P_\rho(b)  +\epsilon}<\infty.\] 
Hence $\{x_w^{1/p}\}\in \ell^p$ and  $\{y_w^{1/q}\}\in \ell^q$  with $p=\frac{1}{t}$ and  $q=\frac{1}{1-t}$  and by Hölder's inequality:
\begin{align*}
( \sum_{w\in G} x_w)^{1/p} .( \sum_{w\in G} y_w)^{1/q} &\geq
 \sum_{w\in G} ((q-1)^{ a \vert w\vert} \rho(w)^{ P_\rho(a)  +\epsilon})^{t} .( (q-1)^{ b \vert w\vert} \rho(w)^{ P_\rho(b)  +\epsilon})^{1-t},\\  
 &\geq \sum_{w\in G} (q-1)^{( t a + (1-t)b) \vert w\vert} \rho(w)^{ (t P_\rho(a)     +(1-t) P_\rho(b)) + \epsilon}.
 \end{align*}  
  
 Finally, by convexity, for all   $t<0$ one has $P_\rho (ts_0) \geq t P_\rho ( s_0) +(1-t) P_\rho ( 0)  .$ In particular, for   $t=1/s_0$ one has
 \[Q(q) \geq Q(2)\cdot (1 -\frac{1}{s_0}) =Q(2) \cdot \big(1+\frac{\log(q-1)}{\tau(2)}\big) .\]

   \bibliography{Biblio}

\begin{thebibliography}{10}

\bibitem{BonkKleinerConfDimGromHypergrps}
M.~Bonk and B.~Kleiner.
\newblock Conformal dimension and {G}romov hyperbolic groups with 2-sphere
  boundary.
\newblock {\em Geom. Topol.}, 9:219--246, 2005.

\bibitem{BourdonImHyperDimConfRigi}
M.~Bourdon.
\newblock Immeubles hyperboliques, dimension conforme et rigidit\'e de
  {M}ostow.
\newblock {\em Geom. Funct. Anal.}, 7(2):245--268, 1997.

\bibitem{BourdonKleinerCLP}
M.~Bourdon and B.~Kleiner.
\newblock Combinatorial modulus, the combinatorial {L}oewner property, and
  {C}oxeter groups.
\newblock {\em Groups Geom. Dyn.}, 7(1):39--107, 2013.

\bibitem{CarrascoConfDimModComb}
M.~Carrasco.
\newblock Conformal dimension and combinatorial modulus of compact metric
  spaces.
\newblock {\em C. R. Acad. Sci. Paris}, Serie I 350:141--145, 2012.

\bibitem{ClaCLPonBRAHB}
A.~Clais.
\newblock Combinatorial {M}odulus on {B}oundary of {R}ight-{A}ngled
  {H}yperbolic {B}uildings.
\newblock {\em Anal. Geom. Metr. Spaces}, 4:Art. 1, 2016.

\bibitem{DavisBook}
M.~W. Davis.
\newblock {\em The geometry and topology of {C}oxeter groups}, volume~32 of
  {\em London Mathematical Society Monographs Series}.
\newblock Princeton University Press, Princeton, NJ, 2008.

\bibitem{DavisMeierTopo}
M.~W. Davis and J.~Meier.
\newblock The topology at infinity of {C}oxeter groups and buildings.
\newblock {\em Comment. Math. Helv.}, 77(4):746--766, 2002.

\bibitem{HaissinskyGeomQConf}
P.~Ha{\"{\i}}ssinsky.
\newblock G\'eom\'etrie quasiconforme, analyse au bord des espaces m\'etriques
  hyperboliques et rigidit\'es [d'apr\`es {M}ostow, {P}ansu, {B}ourdon,
  {P}ajot, {B}onk, {K}leiner{$\ldots$}].
\newblock {\em Ast\'erisque}, (326):Exp. No. 993, ix, 321--362 (2010), 2009.
\newblock S{\'e}minaire Bourbaki. Vol. 2007/2008.

\bibitem{KeithKleiner}
S.~Keith and B.~Kleiner.
\newblock In preparation.

\bibitem{KleinerAsymptoticGeom}
B.~Kleiner.
\newblock The asymptotic geometry of negatively curved spaces: uniformization,
  geometrization and rigidity.
\newblock In {\em International {C}ongress of {M}athematicians. {V}ol. {II}},
  pages 743--768. Eur. Math. Soc., Z\"urich, 2006.

\bibitem{MackayTysonConfDim}
J.~M. Mackay and J.~T. Tyson.
\newblock {\em Conformal dimension}, volume~54 of {\em University Lecture
  Series}.
\newblock American Mathematical Society, Providence, RI, 2010.
\newblock Theory and application.

\bibitem{MeierWhen}
J.~Meier.
\newblock When is the graph product of hyperbolic groups hyperbolic?
\newblock {\em Geom. Dedicata}, 61(1):29--41, 1996.

\bibitem{PansuDimconf}
P.~Pansu.
\newblock Dimension conforme et sph\`ere \`a l'infini des vari\'et\'es \`a
  courbure n\'egative.
\newblock {\em Ann. Acad. Sci. Fenn. Ser. A I Math.}, 14(2):177--212, 1989.

\end{thebibliography}
  \bibliographystyle{abbrv}

\end{document}